\documentclass[a4paper,10pt]{article}
\usepackage{french}
\usepackage[latin1]{inputenc}
\usepackage{amssymb}
\usepackage{amsmath}
\usepackage{diagrams}

\newtheorem{theo}{Théorème}[section]
\newtheorem{lem}[theo]{Lemme}
\newtheorem{prop}[theo]{Proposition}
\newtheorem{cor}[theo]{Corollaire}
\newtheorem{defi}[theo]{Définition}
\newtheorem{rem}[theo]{Remarque}
\newtheorem{nota}[theo]{Notation}
\newtheorem{def-prop}[theo]{Définition-proposition}

\begin{document}
\begin{center}
\bfseries \LARGE La suite spectrale de Cibils pour la
cohomologie des algèbres triangulaires
\end{center}
\vspace{0.5cm}
\begin{center}
\textbf{{\Large Sophie Dourlens}}\\
Laboratoire J.A. Dieudonné\\
Université de Nice Sophia-Antipolis\\
Parc Valrose\\
06108 Nice Cedex 2\\
dourlens@math.unice.fr
\end{center}
\vspace{0.5cm}

\begin{abstract}
Pour étudier la cohomologie de Hochschild d'algèbres
triangulaires $\mathcal T$, 
nous construisons une suite spectrale, dont
les termes sont paramétrés par la longueur des trajectoires 
du carquois associé à $\mathcal T$, et qui converge vers 
$HH^*(\mathcal T)$. Nous en explicitons les composantes, et 
les différentielles au premier niveau qui s'écrivent comme des 
sommes de produits cup. Dans le cas $n=3$, nous étudions les
propriétés de la différentielle au niveau 2.
Enfin, nous appliquons ces résultats à l'algèbre des chemins
d'un carquois sans cycles orientés et les relions à des résultats antérieurs sur
l'algèbre d'incidence d'un complexe simplicial, et plus 
généralement sur l'algèbre des morphismes de certaines 
catégories.
\end{abstract}


\section{Introduction}

D'après Hochschild-Konstant-Rosenberg, si $A$ est l'algèbre
des fonctions régulières sur une variété affine $V$ (en
caractéristique $0$), la cohomologie de
Hochschild de $A$ calcule les sections globales du faisceau des
champs de multi-vecteurs sur $V$. Plus généralement, les
invariants (co)homologiques de l'algèbre $A$ permettent de 
décrire de nombreuses
propriétés géométriques de $V$. 
L'idée de la géométrie non commutative est d'étendre ce type
de résultats au cadre d'algèbres non commutatives. On pourra
se référer, par exemple, à l'introduction de l'article de 
Nest et Tsygan \cite{NT} pour un aperçu synthétique sur ces
questions. Un cas particulièrement intéressant d'algèbres
non commutatives est celui des algèbres triangulaires, qui
jouent, dans ce contexte (au moins si leur sous-algèbre
diagonale est commutative), le r\^ole que jouent les algèbres de Lie
nilpotentes en théorie des algèbres de Lie. Comme nous le
verrons, les propriétés de nilpotence des algèbres de Lie
triangulaires se reflètent d'ailleurs dans la suite spectrale de la
cohomologie de Hochschild d'algèbres triangulaires,
introduite et \'etudi\'ee dans cet article.\\

Ces groupes de cohomologie ont été très étudiés
récemment dans le cas d'algèbres de taille 2 \cite{Cib, 
CMRS, GG, MP}. Rappelons en particulier les résultats de
Cibils \cite{Cib}: la cohomologie de l'algèbre $\mathcal{T}=\left(\begin{array}{cc}
A_{1} &\\M_{1} & A_{2} \end{array}\right)$ est reli\'ee 
 \`a la cohomologie de 
Hochschild de $A=A_{1}\times A_{2}$ et aux
groupes d'extension de $M$ par une suite exacte
longue de cohomologie. De plus, le morphisme de
connexion s'exprime comme un produit \emph{cup}.\\

Nous généralisons ici la méthode de Cibils aux alg\`ebres
triangulaires de tailles sup\'erieures. Rappelons tout d'abord 
la construction d'une telle algèbre. 
Soient $A_1$,\ldots, $A_n$ des algèbres 
unitaires sur un corps $k$, et soit, pour tous $1\leq i<j\leq n$, un 
$A_j$-$A_i$-bimodule non nul noté $_jM_i$. On suppose que la collection 
de bimodules $\{_jM_i\}$ est munie de morphismes 
$\mu_{l,j,i}:\ {_lM_j}\otimes{_jM_i}\longrightarrow
{_lM_i}$ pour tous $1\leq i<j<l\leq n$ (le symbole $\otimes$
d\'esignant le produit tensoriel sous $k$), 
tels que les diagrammes
\begin{diagram}
{_mM_l}\otimes{_lM_j}\otimes{_jM_i}&\rTo^{1_{_mM_l}\otimes\mu_{l,j,i}}&{_mM_l}\otimes{_lM_i}\\
\dTo^{\mu_{m,l,j}\otimes 1_{_jM_i}}&\circlearrowleft&\dTo_{\mu_{m,l,i}}\\
{_mM_j}\otimes{_jM_i}&\rTo_{\mu_{m,j,i}}&{_mM_i}
\end{diagram}
soient commutatifs, quels que soient $1\leq i<j<l<m\leq n$. 
L'\emph{algèbre triangulaire} $\mathcal T$ associée à ces 
bimodules est
$$\mathcal T\qquad=\qquad\left(\begin{array}{cccccc}
A_1\\
_2M_1&A_2\\
_3M_1&_3M_2&A_3\\
\\
\vdots&\vdots&&\ddots&\ddots\\
\\
_nM_1&_nM_2&\cdots&\cdots&_nM_{n-1}&A_n
\end{array}\right),$$
la somme et le produit de deux éléments étant donnés par
la somme et le produit des matrices carrées usuelles. 
La structure multiplicative de $\mathcal T$ est donn\'ee par
les morphismes $\mu_{l,j,i}$, et la commutativité
des diagrammes ci-dessus garantit l'associativité du produit 
de $\mathcal T$. 
Un cas particulier d'algèbres triangulaires est celui des
\emph{algèbres triangulaires tensorielles}: les modules
${_{i+\alpha}M_i}$ pour $\alpha\geq 2$ sont alors les produits
tensoriels ${_{i+\alpha}M_{i+\alpha-1}}\otimes_{A_{i+\alpha-1}}
\ldots\otimes_{A_{i+1}}{_{i+1}M_i}$. 
Le r\^ole des morphismes $\mu_{l,j,i}$ est tenu
dans ce cas par les surjections canoniques 
${_lM_j}\otimes_k\,{_jM_i}\twoheadrightarrow
{_lM_j}\otimes_{A_j}{_jM_i}$, et l'associativité du
produit de $\mathcal T$ r\'esulte de l'associativité du
produit tensoriel.\\

Suivant \cite{Cib}, pour étudier la cohomologie d'une algèbre
triangulaire $\mathcal T$ de taille $n$, nous utiliserons son 
\emph{carquois}, qui est de la forme 
$$\mathcal Q\ :\qquad e_1\longrightarrow e_2\longrightarrow\cdots
\longrightarrow e_n\ .$$
Les sommets de $\mathcal Q$ correspondent aux 
$k$-algèbres $A_i$, les flèches aux bimodules 
$_{i+1}M_i$, et les chemins de longueur $\alpha$
($\alpha\geq 2$) 
aux bimodules $_{i+\alpha}M_i$. 
La cohomologie de Hochschild de $\mathcal T$ s'exprime alors
à l'aide d'un complexe de Hochschild relatif \cite{GS}, dont les 
composantes sont paramétrées
par les \emph{trajectoires} du carquois $\mathcal Q$ de
$\mathcal T$. Cette paramétrisation
induit une filtration du complexe par la 
\emph{longueur} des trajectoires de $\mathcal Q$, et donc une 
suite spectrale, qui  converge
vers la cohomologie de Hochschild de
$\mathcal{T}$. 
Le propos de cet article est d'étudier et de
décrire cette suite spectrale (au-delà du cas $2\times 2$
qui, on l'a dit, a déjà été amplement étudié). Comme on
pouvait l'espérer, les diff\'erentielles au niveau 1 
s'expriment à nouveau, dans le cas général, comme des produits \emph{cup}. Dans le
cas des algèbres triangulaires tensorielles de taille 3, 
nous montrons que
la différentielle au niveau 2 a des composantes nulles;
nous décrivons un cas particulier dans lequel la suite
spectrale dégénère au niveau 2.\\

Nous appliquons enfin ces résultats à l'algèbre des chemins
d'un carquois sans cycles orientés, qui est triangulaire.
Nous explicitons sur un exemple de carquois le calcul de la
suite spectrale. Si le carquois provient d'un complexe
simplicial, nous retrouvons un résultat de Cibils
\cite{Cib2}. D'autre part, en considérant un carquois comme
une catégorie, nous complétons des résultats de Bendiffalah
et Guin \cite{BG} sur la cohomologie de certaines catégories.

\section{Cohomologie de Hochschild relative}\label{relative}

La clé du calcul de la cohomologie de Hochschild d'une
algèbre triangulaire est la dé\-ter\-mi\-na\-tion d'un complexe
plus petit que le complexe de Hochschild classique, dont la
cohomologie est la cohomologie de Hochschild de l'algèbre
triangulaire. Nous
donnons dans cette section la construction d'un tel
complexe, qui s'interprète comme un complexe de 
Hochschild relatif \cite{GS}.

Soit $\mathcal T$ l'algèbre triangulaire définie dans
l'introduction. 
Les algèbres $A_i$ étant supposées unitaires, $\mathcal T$ a une
sous-algèbre remarquable
$$R\quad=\quad\left(\begin{array}{ccc}
k\\
&\ddots\\
&&k
\end{array}\right).$$
Pour $i=1\ldots,n$, on note $e_i\in R$ l'image de $1_{A_i}$ 
par l'injection $A_i\hookrightarrow\mathcal T$. Ces 
éléments forment un système $\{e_1,\ldots,e_n\}$ d'idempotents 
orthogonaux, tel que $\sum_ie_i=1_\mathcal T$. De plus, ce
système est une $k$-base de $R$: tout élément
$r$ de $R$ s'écrit de manière unique $r=r\Big(\sum_ie_i\Big)
=\sum_ire_i=\sum_ir_ie_i$, où les 
$r_i$ sont des éléments de $k$.
Soit $R^e=R\otimes_kR^{op}$ l'algèbre enveloppante de $R$.
Cette alg\`ebre possède un élément remarquable
$e$, d\'efini par:
$$e:=\sum_{i=1}^ne_i\otimes e_i^{op}.$$

\begin{lem}\label{adm}
Pour tout morphisme de $\mathcal T$-bimodules
$f:X\longrightarrow Y$, il existe $\sigma:
Y\longrightarrow X$ $R$-linéaire tel que $f\sigma f=f$.
\end{lem}

PREUVE. Soit la décomposition canonique de $f$ en une
surjection suivie d'une injection $f:X\stackrel{p}
{\twoheadrightarrow}Im\,f\stackrel{i}{\hookrightarrow}Y$. 
Comme $k$ est un corps, la surjection $p$ admet une section
$k$-linéaire, et l'injection $i$ admet une rétraction
$k$-linéaire. La composée de ces deux applications donne une
application $k$-linéaire $\sigma:Y\longrightarrow X$ telle que 
$f\sigma f=f$.
Nous allons construire à partir de $\sigma$ une
application $\sigma'$ $R$-linéaire ayant la m\^eme propriété. 
Pour cela, notons $S=R^e$, et considérons
l'automorphisme de $S^e$ suivant:
\begin{eqnarray*}
S^e=R\otimes R^{op}\otimes R^{op}\otimes R
&\stackrel {\cong}{\longrightarrow}& R\otimes R^{op}\otimes
R^{op}\otimes R=S^e\\
p\otimes q^{op}\otimes r^{op}\otimes s &\longmapsto& p\otimes 
r^{op}\otimes q^{op}\otimes s.
\end{eqnarray*}
Soit $e_S$ l'image de $e\otimes e^{op}$ par cet automorphisme.
Précisément, $e_S$ s'écrit
$$e_S=\sum_{i,j=1}^n e_i\otimes e_j^{op}\otimes e_i^{op}\otimes
e_j.$$
Remarquons que 
$Hom_k(Y,X)$ a une stucture de $S^e$-module à gauche: 
si $\sigma\in Hom_k(Y,X)$,
et $p, q, r, s\in R$, alors 
l'action de $(p\otimes q^{op}\otimes r^{op}\otimes s)$ sur 
$\sigma$ est donnée par
\begin{eqnarray*}
(p\otimes q^{op}\otimes r^{op}\otimes s).\sigma:
Y&\longrightarrow&X\\
y&\longmapsto&p\big(\sigma(rys)\big)q.
\end{eqnarray*}
Posons alors $\sigma'=e_S.\sigma$. D'une part, $\sigma'$ est un
morphisme de $R$-bimodules: en effet, pour $r, s \in R$, et $y\in
Y$, la décomposition de $r$ et $s$ dans la base
$\{e_1,\ldots,e_n\}$ et la $k$-linéarité de $\sigma$ permettent
d'écrire
\begin{eqnarray*}
\sigma'(rys)&=&\sum_{i,j=1}^ne_i\sigma(e_iryse_j)e_j
=\sum_{i,j=1}^ne_i\sigma(r_ie_iye_js_j)e_j
=\sum_{i,j=1}^ne_ir_i\sigma(e_iye_j)s_je_j\\
&=&\sum_{i,j=1}^nre_i\sigma(e_iye_j)e_js
=r\Big(\sum_{i,j=1}^ne_i\sigma(e_iye_j)e_j\Big)s
=r\sigma'(y)s.
\end{eqnarray*}
D'autre part, $\sigma'$ vérifie la propriété demandée, à savoir
$f\sigma'f=f$: en effet, comme $f$ est un morphisme de
$\mathcal T$-bimodules, il commute en particulier avec les
$e_i$; on a ainsi pour $x\in X$
\begin{eqnarray*}
f\sigma'f(x)&=&f\Big(\sum_{i,j=1}^ne_i\sigma(e_if(x)e_j)e_j\Big)
=\sum_{i,j=1}^ne_if(\sigma(f(e_ixe_j)))e_j\\
&=&\sum_{i,j=1}^ne_if(e_ixe_j)e_j
=\sum_{i,j=1}^ne_if(x)e_j=\Big(\sum_{i=1}^ne_i\Big)f(x)\Big(\sum_{j=1}^ne_j\Big)=f(x).
\ \diamondsuit
\end{eqnarray*}

\begin{prop}\label{resolution}
Le complexe suivant
$$ \ldots \rightarrow \mathcal T^{\otimes_{R}l+2} \stackrel
{d_{l}}
\longrightarrow \mathcal T^{\otimes_{R}l+1}\longrightarrow
\ldots \longrightarrow \mathcal T\otimes_{R}\mathcal T\longrightarrow \mathcal T
\rightarrow 0$$
avec la diff\'erentielle
$$d_{l}(t_{1}\otimes\ldots\otimes t_{l+2})=
\sum_{i=1}^{l+1}(-1)^{i}~
t_{1}\otimes\ldots\otimes
t_{i}t_{i+1}\otimes\ldots\otimes t_{l+2}$$
est une r\'esolution projective de $\mathcal T$ par des
$\mathcal T$-bimodules.
\end{prop}

PREUVE.
Tout d'abord, on a une homotopie
\begin{eqnarray*}
s_{l}: & T^{\otimes_{R}l+1} & \longrightarrow \quad
T^{\otimes_{R}l+2}\\
& t_{1}\otimes\ldots\otimes t_{l+1} & \longmapsto
\quad 1\otimes t_{1}\otimes\ldots\otimes t_{l+1},
\end{eqnarray*}
qui donne l'exactitude du complexe. Il reste donc à montrer que
$\mathcal T^{\otimes_Rl+2}$ est projectif pour $l\geq 0$.
Soient, dans la catégorie des $\mathcal T$-bimodules, un
morphisme surjectif
$\pi:X\twoheadrightarrow Y$ et un morphisme $f:\mathcal T
^{\otimes_Rl+2} \rightarrow Y$. Il s'agit de montrer qu'il
existe $F:\mathcal T^{\otimes_Rl+2} \rightarrow X$ faisant
commuter le diagramme
\begin{diagram}
& & \mathcal T^{\otimes_Rl+2}&\\
& \ldDotsto^{\exists ?\ F}_{\circlearrowleft} & \dTo^{f}\\
X & \rTo^{\pi} & Y & \rTo & 0\\
\end{diagram}
Or, d'après  la relation d'adjonction
$$Hom_{\mathcal T-\mathcal T}(\mathcal T^{\otimes_Rl+2},Y)\cong
Hom_{R-R}(\mathcal T^{\otimes_Rl},Y),$$
le diagramme précédent induit dans la catégorie des
$R$-bimodules un nouveau diagramme:
\begin{diagram}
& & \mathcal T^{\otimes_Rl}&\\
&  & \dTo^{\overline f}\\
X & \rTo^{\pi} & Y & \rTo & 0\\
\end{diagram}
D'après le lemme \ref{adm}, il existe un morphisme de
$R$-bimodules $\sigma:Y\rightarrow X$ tel que
$\pi\sigma\pi=\pi$. Comme $\pi$ est surjectif, cette égalité
devient $\pi\sigma=id_Y$. Le morphisme $\overline
F=\sigma\overline f$ fait donc commuter le diagramme
\begin{diagram}
& & \mathcal T^{\otimes_Rl}&\\
& \ldDotsto^{\overline F}_{\circlearrowleft} & \dTo^{\overline f}\\
X & \pile{\rTo^{\pi}\\ \lDotsto_{\sigma}} & Y & \rTo & 0\\
\end{diagram}
En appliquant à nouveau la relation d'adjonction, on obtient le morphisme
de $\mathcal T$-bimodules $F:\mathcal
T^{\otimes_Rl+2}\longrightarrow X$ qui donne la projectivité de
$\mathcal T^{\otimes_Rl+2}$. $\diamondsuit$

\begin{cor}\label{complexe}
Le complexe suivant
$$\mathcal{H}:\qquad 0\rightarrow X^{R}\longrightarrow
Hom_{R-R}(\mathcal T,X)\longrightarrow
\ldots\longrightarrow Hom_{R-R}\big(
\mathcal T^{\otimes_{R}l},X\big)
\stackrel{\delta _{l}}\longrightarrow \ldots$$
avec le cobord:
$$\begin{array}{l}
\ \delta _{1}x(t)=tx-xt\\
\\
\begin{array}{rcl}
\delta _{l}f(t_{1}\otimes\ldots\otimes t_{l+1}) & = &
t_{1}~f(t_{2}\otimes\ldots\otimes t_{l+1})\\
&&+ \sum\limits_{i=1}^{l}(-1)^{i}~f(t_{1}\otimes\ldots\otimes
t_{i}t_{i+1}\otimes\ldots\otimes t_{l+1})\\
&&+ (-1)^{l+1}f(t_{1}\otimes\ldots\otimes t_{l})t_{l+1}.
\end{array}
\end{array}$$
calcule $HH^{*}(\mathcal T,X)$.
\end{cor}

PREUVE.
Il suffit d'appliquer le foncteur $Hom_{\mathcal T-\mathcal
T}(\ -\ ,X)$ à la résolution projective de $\mathcal T$ donnée
par la proposition précédente, puis d'utiliser les relations
$$\left\{ \begin{array}{l}
Hom_{\mathcal T-\mathcal T}\big(\mathcal T^{\otimes_{R}l+2},\ X\big) \cong 
Hom_{R-R}\big(\mathcal T^{\otimes_{R}l},\ X\big)\quad
\forall\,l>0\,,\\
Hom_{\mathcal T-\mathcal T}\big(\mathcal T\otimes_{R}\mathcal T,\ X\big) \cong
X^{R}.\ \diamondsuit
\end{array}\right.$$
\begin{rem}\emph{
Les résultats de cette section sont un cas particulier
d'une théorie de cohomologie relative développée par 
Gerstenhaber et Schack \cite{GS} (section 1). Nous avons en
fait construit le complexe de Hochschild de $\mathcal T$ 
\emph{relatif} à la sous-algèbre $R$. La résolution donnée à la
proposition \ref{resolution} est une résolution
projective de $\mathcal T$ \emph{relativement à $R$}, ou
résolution \emph{$R$-projective} de $\mathcal T$. Les
propriétés de l'élément $e\in R^e$ font de $R$ une algèbre
\emph{séparable}; c'est pourquoi une résolution
$R$-projective de $\mathcal T$ est en fait une résolution 
projective de $\mathcal T$ (lemme \ref{adm}). Ainsi, la 
cohomologie du complexe relatif co\"\i ncide avec la 
cohomologie de Hochschild de $\mathcal T$. Cibils a donn\'e
dans \cite{Cib} une autre d\'emonstration de la proposition
\ref{resolution}, sans utiliser la notion de cohomologie 
relative de Gerstenhaber et Schack.}\end{rem}

\section{Le carquois d'une algèbre triangulaire}\label{carquois}

Le but de cette section est de fournir une paramétrisation des
termes de la résolution projective donnée dans la section
précédente (propostion \ref{resolution}). Pour cela, nous
allons rappeler la notion de carquois dans le cas
d'une algèbre triangulaire tensorielle \cite{Cib}, puis
l'adapter dans le cas d'une algèbre triangulaire quelconque.

\begin{defi}
Soit $A$ un anneau. Un \emph{syst\`eme $Q_{0}$
d'idempotents centraux de $A$} est un ensemble fini
d'idempotents non nuls, centraux, orthogonaux, tel que
\mbox{$\sum_{e \in Q_{0}} e = 1 $.} 
Le \emph{carquois} d'un $A$-bimodule $M$ non nul
associ\'e \`a $Q_{0}$ est un graphe orient\'e dont
les sommets sonts les \'el\'ements de $Q_{0}$; il
existe une fl\`eche entre le sommet $e$ et le sommet
$e'$ si et seulement si $e'Me\neq0$.
\end{defi}

Soit $\mathcal T$ l'algèbre triangulaire \emph{tensorielle} de taille $n$ définie
dans l'introduction. On pose $M={_2M_1}\oplus\ldots\oplus 
{_nM_{n-1}}$, et $A=A_{1}\times\ldots\times A_{n}$. $M$ poss\`ede
naturellement une structure de $A$-bimodule, induite
par les actions de $A_{i+1}$ et $A_{i}$ sur chaque
$_{i+1}M_{i}$, et l'action nulle pour les autres facteurs de
$A$. 
$A$ poss\`ede un syst\`eme naturel d'idempotents 
centraux $Q_{0}=\{e_{1},\ldots,e_{n}\}$, les $e_i$ ayant été
définis précédemment comme les images de $1_{A_i}$ par
l'injection $A_i\hookrightarrow A$.  
Le carquois de $M$ par rapport \`a $Q_{0}$ est
$$e_{1}\longrightarrow e_{2}\longrightarrow\cdots
\longrightarrow e_{n}.$$
Nous ferons référence à ce carquois comme le \emph{carquois de
l'algèbre $\mathcal T$}. 
Les $k$-algèbres $A_i$
correspondent aux sommets $e_i$ du carquois, 
les bimodules $_{i+1}M_i$ correspondent aux flèches, et les produits
tensoriels de la forme $_{k+\alpha}M_{k+\alpha-1}\otimes_
{A_{k+\alpha-1}}\ldots\otimes_{A_{k+1}}{_{k+1}M_k}$
correspondent aux chemins de longueur $\alpha$ 
dans le carquois.

\begin{defi}
Une \emph{$l$-trajectoire} d'un carquois $\mathcal{Q}$ est une 
suite $(\gamma_{l}, \ldots, \gamma_{1})$ de l chemins
cons\'ecutifs de $\mathcal{Q}$, sans restriction sur la
longueur des chemins (les sommets du carquois, qui sont des
chemins de longueur nulle, peuvent \^etre inclus dans la composition
d'une trajectoire). 
On appelle \emph{longueur d'une $l$-trajectoire} le
nombre de chemins \emph{de longueur non nulle} qui la
composent. On note $|\tau|$ la longueur de la $l$-trajectoire $\tau$; on a
$0\leq |\tau|\leq l$.
\end{defi}

L'ensemble des $l$-trajectoires est noté
$TR_{l}(\mathcal{Q})$. Nous reprenons les notations de Cibils
\cite{Cib} associées au carquois $\mathcal{Q}$ d'un $A$-
bimodule $M$:
\begin{itemize}
\item pour chaque sommet $e$ de $\mathcal{Q}$, on pose 
$M_{e}=Ae$ (c'est une sous-algèbre de $A$);
\item si $a$ est une fl\`eche de $\mathcal{Q}$, on pose
$M_{a}=t(a)Ms(a)$, o\`u $t(a)$ et $s(a)$ d\'esignent
respectivement le but et la source de la fl\`eche
$a$; $M_{a}$ est non nul par d\'efinition des
fl\`eches de $\mathcal{Q}$;
\item si $\gamma = a_{p}\ldots a_{1}$ est un chemin de
$\mathcal{Q}$, alors on pose $M_{\gamma} = M_{a_{p}}\otimes_
{As(a_{p})}\ldots\otimes_{At(a_{1})} M_{a_{1}}$;
\item enfin, si $\tau = (\gamma_{l},\ldots,\gamma_{1})$
est une $l$-trajectoire, alors on pose $M_{\tau} =
M_{\gamma_{l}}\otimes_{k}\ldots\otimes_{k}
M_{\gamma_{1}}.$\\
\end{itemize}

Soient $A$ une $k$-algèbre munie d'un système d'idempotents
centraux $Q_0$, et $M$ un $A$-bimodule de carquois $\mathcal
Q$ associé à $Q_0$. D'après un résultat de Cibils \cite{Cib}, 
si $R=\prod_ik\,e_i$ est la sous-$k$-algèbre de $A$ engendrée par 
$Q_0$, et si $T_A(M)$ est
l'algèbre tensorielle de $M$ sous $A$, alors la $l^e$ puissance
tensorielle de $T_A(M)$ sous $R$ se décompose selon les l-trajectoires
de $\mathcal Q$:
$$\Big(T_A(M)\Big)^{\otimes_{R}l} = \bigoplus_{\tau\in
TR_{l}(\mathcal{Q})}
M_{\tau}\ .$$
La proposition suivante permet d'appliquer ce résultat à
l'algèbre triangulaire tensorielle $\mathcal T$.

\begin{prop}\label{iso}
Si $A=A_1\times\ldots\times A_n$, et
$M={_2M_1}\oplus\ldots\oplus {_nM_{n-1}}$, alors les alg\`ebres 
$\mathcal{T}$ et $T_{A}(M)$ sont isomorphes.
\end{prop}

PREUVE.
On a, par d\'efinition de l'action de $A$ sur $M$:
$$M\otimes_{A}M=\bigg(\bigoplus_{i=1}^{n-1}
{_{i+1}M_{i}}\bigg)\otimes_{A}\bigg(\bigoplus_{i=1}^{n-1}{_{i+1}M_{i}}\bigg) 
=\bigoplus_{i,j=1}^{n-1}{_{j+1}M_j}\otimes_A{_{i+1}M_i},$$
et
$$_{j+1}M_j\otimes_{A}{_{i+1}M_i}=\left\{ \begin{array}{ll}
0 & \textrm{si}\ j\neq i+1,\\
{_{i+2}M_{i+1}}\otimes_{A_{i+1}}{_{i+1}M_i} & \textrm{sinon}.
\end{array} \right.$$
On en d\'eduit:
$$M^{\otimes_{A}r}=\bigoplus_{i=1}^{n-r-1}
{_{i+r+1}M_{i+r}}\otimes_{A_{i+r}}{_{i+r}M_{i+r-1}}\otimes_{A_{i+r-1}}
\ldots\otimes_{A_{i+1}}{_{i+1}M_i}\quad \text{pour } 1\leq r\leq
n-1,$$
et $M^{\otimes_{A}n}=0$, ce qui donne l'isomorphisme
annoncé.$\ \diamondsuit$\\

Dans le cas d'une algèbre triangulaire quelconque $\mathcal
T$ de taille $n$, par analogie avec 
le cas tensoriel, on peut aussi lui associer le carquois
$$\mathcal Q\ :\qquad e_1\longrightarrow e_{2}\longrightarrow\cdots
\longrightarrow e_{n}.$$
Le sommet $e_i$ de $\mathcal Q$ correspond toujours à l'algèbre
$A_i$, la flèche entre le sommet $e_i$ et $e_{i+1}$ correspond 
au bimodule $_{i+1}M_i$, et le chemin $\gamma$ de longueur 
$\alpha\geq 2$ entre les sommets $e_i$ et $e_{i+\alpha}$
correspond maintenant au bimodule $M_\gamma={_{i+\alpha}M_{i}}$.  
Si $\tau=(\gamma_l,\ldots,\gamma_1)$ est une $l$-trajectoire,
le module correspondant est toujours
$M_\tau=M_{\gamma_l}\otimes_k\ldots\otimes_kM_{\gamma_1}$. 
Ainsi, la $l^e$ puissance
tensorielle de $\mathcal T$ sous $R$ peut encore \^etre
exprimée en fonction des $l$-trajectoires de $\mathcal Q$:
$$\mathcal T^{\otimes_{R}l} = \bigoplus_{\tau\in
TR_{l}(\mathcal{Q})}
M_{\tau}\ .$$

Cette décomposition de $\mathcal T^{\otimes_{R}l}$ par les
$l$-trajectoires de $\mathcal Q$ nous donne
donc une pa\-ra\-mé\-tri\-sa\-tion des termes de la résolution
projective de $\mathcal T$ donnée par la proposition
\ref{resolution}, $R=\prod_ik\,e_i$ étant bien la sous-algèbre séparable de 
$\mathcal T$ utilisée dans la section précédente:
$$R\quad=\quad\left(\begin{array}{ccc}
k\\
&\ddots\\
&&k
\end{array}\right).$$

\section{Une suite spectrale pour la cohomologie des algèbres
triangulaires}\label{suite spectrale}

\subsection{La filtration du complexe de Hochschild relatif}

Considérons le complexe $\mathcal H$ défini dans le
corollaire \ref{complexe}, dont la cohomologie 
est $HH^*(\mathcal T,X)$:
$$\mathcal{H}:\qquad 0\rightarrow X^{R}\longrightarrow
Hom_{R-R}(T,X)\longrightarrow\ldots
\longrightarrow Hom_{R-R}\big(
T^{\otimes_{R}l},X\big)
\stackrel{\delta _{l}}\longrightarrow \ldots$$
D'après la décomposition de $\mathcal T^{\otimes_{R}l}$ en
fonction des $l$-trajectoires du carquois $\mathcal Q$ de
$\mathcal T$, le $l^e$ terme de $\mathcal H$ s'écrit
$$\mathcal H^l=Hom_{R-R}\Bigg(\bigoplus_{\tau\in TR_{l}
(\mathcal{Q})}M_{\tau}\ ,X\Bigg).$$
On en d\'eduit une filtration du complexe $\mathcal{H}$ param\'etr\'ee
par la longueur des trajectoires de $\mathcal Q$; on a pr\'ecis\'ement
$$0=F^{n}\mathcal{H}^{l}\subset F^{n-1}\mathcal{H}^{l}\subset\ldots\subset
F^{t}\mathcal{H}^{l}\subset\ldots\subset F^{1}\mathcal{H}^{l}\subset
F^{0}\mathcal{H}^{l}=\mathcal{H}^{l},$$
avec
$$F^{t}\mathcal{H}^{l}=Hom_{R-R}\Bigg( \bigoplus_{\substack{\tau\in
TR_{l}(\mathcal{Q}) \\ |\tau |\geq t}} M_{\tau}\ ,\ X\Bigg).$$

La notation suivante permet d'exprimer le module associé à une
trajectoire en fonction des composantes de l'algèbre 
triangulaire:

\begin{nota}
Soient des entiers ordonnés $1\leq k_1<k_2<\ldots<k_{t+1}\leq
n$, et des entiers positifs ou nuls $p_1,\ldots,p_{t+1}$. 
On note $\ ^{p_{t+1}}{_{k_{t+1}}M_{k_t}}^{p_t}\ldots
{^{p_2}{_{k_2}M_{k_1}}^{p_1}}$
le produit tensoriel sous le corps $k$ : 
$A_{k_{t+1}}^{\otimes p_{t+1}}\otimes\,{_{k_{t+1}}M_{k_t}}^{p_t}\otimes
A_{k_t}^{\otimes p_t}\otimes\ldots\otimes A_{k_2}^{\otimes
p_2}\otimes\,{_{k_2}M_{k_1}}\otimes A_{k_1}^{\otimes p_1}$.\\
Pour un $\mathcal T$-bimodule $X$, on note
$_{j}X_{i}:=e_{j}Xe_{i}$, et
$X_{0}:=\bigoplus\limits_{1\leq i\leq
n}{_{i}X_{i}}.$
\end{nota}

Nous allons étudier la suite spectrale associ\'ee \`a
la filtration ci-dessus; elle converge \textit{a priori} au niveau 
$n$ vers $HH^{*}(\mathcal{T},X).$ 
Afin de faciliter la lecture, nous allons traiter le cas
$n=3$ en d\'etail. Dans un second temps, nous
formulerons les r\'esultats pour $n$ quelconque.

\subsection{La suite spectrale pour $n=3$}

Nous considérons donc
$${\mathcal{T}} \quad = \quad
\left(\begin{array}{ccc}
A_1\\
_2M_1 & A_2\\
_3M_1 & _3M_2 & A_3
\end{array}\right),$$
avec une application $k$-linéaire $\mu:{_3M_2}\otimes{_2M_1}
\longrightarrow{_3M_1}$ qui donne la structure
multiplicative de $\mathcal T$. 
Avec la notation ci-dessus, la filtration du complexe 
$\mathcal{H}$ devient
$$\left\{\begin{array}{l}
F^{0}\mathcal{H}^{l}=\mathcal{H}^{l}=Hom_{R-R}(\mathcal{T}^{\otimes_{R}l},X)\\
\\
F^{1}\mathcal{H}^{l}=Hom_{R-R}\big(\bigoplus\limits
_{p+q=l-1}\!\!\!{^q{_2M_1}^p}\oplus\bigoplus\limits_
{p+q=l-1}\!\!\!{^q{_3M_2}^p}\oplus\bigoplus\limits_
{p+q=l-1}\!\!\!\!{^q{_3M_1}^p}\oplus\!\!\bigoplus\limits_
{p+q+r=l-2}\!\!\!\!\!\!\!\!{^{r}{_3M_2}^{q}{_2M_1}^{p}}\,,\
X\ \big)\\
\\
F^{2}\mathcal{H}^{l}=Hom_{R-R}\big(\bigoplus\limits_
{p+q+r=l-2}{^{r}{_3M_2}^{q}{_2M_1}^{p}},X\ \big)\\
\\
F^{t}\mathcal{H}^{l}=0\quad si\ t\geq 3\,.
\end{array}\right.$$
Les termes au niveau 0 de la suite spectrale associée sont
donc
$$\left\{\begin{array}{l}
E_{0}^{0,l}\ =\ Hom_{R-R}\big(A^{\otimes_{R}l}\ ,\ X\big)\
,\quad\textrm{où}\ A=A_1\times\ldots\times A_n\\
\\
E_{0}^{1,l}\ =\ Hom_{R-R}\big(\bigoplus\limits
_{p+q=l}\!\!\!{^{q}{_2M_1}^{p}}\oplus\bigoplus\limits_
{p+q=l}\!\!\!{^{q}{_3M_2}^{p}}\oplus
\bigoplus\limits_{p+q=l}\!\!\!\!
{^{q}{_3M_1}^{p}}\ ,\ X\ \big)\\
\\
E_{0}^{2,l}\ =\ Hom_{R-R}\big(\bigoplus\limits_
{p+q+r=l}{^{r}{_3M_2}^{q}{_2M_1}^{p}}\ ,\ X\ \big)\,.
\end{array}\right.$$
Nous allons maintenant décrire le niveau 1 de la suite
spectrale : nous allons calculer les trois termes, puis 
les deux différentielles. 
Tout d'abord, le terme $E_{1}^{0,*}$
est $HH^*(A,{X_0})=\bigoplus\limits_{i=1}^{n}HH^*(A_i,{_iX_i})$; 
on a en effet
$A^{\otimes_{R}l}=\bigoplus\limits_{i=1}^{n}A_i^{\otimes_kl}$, 
et $Hom_{R-R}\big(A_i^{\otimes_kl}\ ,\ X\big)\cong
Hom_{k}\big(A_i^{\otimes_{k}l} ,\ {_iX_i}\big).$

\begin{prop}\label{prop1}
Le terme $E_{1}^{1,*}$ est
$$Ext^{*}_{A_2-A_1}({_2M_1},{_2X_1})
\oplus Ext^{*}_{A_3-A_2}({_3M_2},{_3X_2})
\oplus Ext^{*}_{A_3-A_1}({_3M_1},{_3X_1}).$$
\end{prop}

La démonstration de cette proposition repose sur le lemme
suivant:

\begin{lem}\label{lem1}\emph{\cite{Cib}}
Soient $B$ et $A$ des $k$-alg\`ebres, et $N$ un
$B$-$A$-bimodule.\\
Le complexe suivant
$$\ldots\rightarrow\bigoplus_{\substack{{p+q=l+2}\\{p>0,\,q>0}}}
{^{q}N^{p}}\stackrel{d_{l}}{\longrightarrow}\bigoplus
_{\substack{{p+q=l+1}\\{p>0,\,q>0}}}{^{q}N^{p}}\rightarrow
\ldots\rightarrow{^{1}N^{1}}\stackrel{\epsilon}
{\longrightarrow} N\rightarrow 0$$
avec la diff\'erentielle
\begin{eqnarray*}
d_{l}:\qquad {^{q}N^{p}} & \longrightarrow &
{^{q-1}N^{p}}\oplus{^{q}N^{p-1}}\\
(b_{1},\ldots,b_{q},x,a_{1},\ldots,a_{p}) & \longmapsto
& \quad \sum_{i=1}^{p-1}(-1)^{i+1}\ (b_{1},\ldots,b_{i}b_{i+1},\ldots,b_{q},x,a_{1},\ldots,a_{p})\\
& & +\ (-1)^{p+1}\ (b_{1},\ldots,b_{q-1},b_{q}x,a_{1},\ldots,a_{p})\\
& & +\ (-1)^{p+2}\ (b_{1},\ldots,b_{q},xa_{1},a_{2},\ldots,a_{p})\\
& & +\ \sum_{i=1}^{q-1}(-1)^{p+i+2}\ (b_{1},\ldots,b_{q},x,a_{1},\ldots,a_{i}a_{i+1},\ldots,a_{p})
\end{eqnarray*}
\begin{tabular}{rll}
$et\qquad \epsilon : \quad {^{1}N^{1}}$ & $\longrightarrow$ & $N$\\
$(b,x,a)$ & $\longmapsto$ & $bxa$
\end{tabular}\\
\\
est une r\'esolution libre de $N$.
\end{lem}

PREUVE DE LA PROPOSITION~\ref{prop1}.
Pour obtenir la première composante de $E_1^{1,l}$, il 
suffit de remplacer dans le lemme~\ref{lem1} 
$N$ par $_2M_1$, puis d'appliquer le foncteur
$Hom_{A_2-A_1}(\,-\,,{_2X_1})$ \`a la r\'esolution
libre de $_2M_1$. Par adjonction, on obtient
$$Hom_{A_2-A_1}\Bigg(\bigoplus_{\substack{
{p+q=l+2}\\{p>0,\,q>0}}}
{^q{_2M_1}^p},{_2X_1}\Bigg)\cong
Hom_k\Bigg(
\bigoplus_{p+q=l}{^q{_2M_1}^p},{_2X_1}\Bigg).$$
Enfin, il reste à remarquer que l'on a aussi un isomorphisme
$$Hom_k\Bigg(
\bigoplus_{p+q=l}{^q{_2M_1}^p},{_2X_1}\Bigg)\cong 
Hom_{R-R}\Bigg(
\bigoplus_{p+q=l}{^q{_2M_1}^p},X\Bigg).$$
On trouve ainsi la première composante de $E_0^{0,l}$. Les
deux autres composantes se traitent de la m\^eme manière. 
$\diamondsuit$

\begin{prop}\label{Ext}
Le terme $E_{1}^{2,*}$ est $\mathbf{Ext}
^{*}_{A_{3}-A_{1}}(\mathcal{C},{_{3}X_{1}})$,
o\`u  $\mathcal{C}$ est un complexe de
cha\^\i nes dont l'homologie est 
$Tor_{*}^{A_{2}}({_3M_2},{_2M_1})$.
\end{prop}

On peut se reporter à \cite{Weibel} pour la notion de
foncteur \emph{HyperExt}, $\mathbf{Ext}^{*}$. En particulier,
$E_1^{2,*}$ est 
l'aboutissement d'une suite spectrale
dont les termes au niveau 2 sont
$$\mathcal E_2^{p,q}=Ext^p_{A_3-A_1}(Tor_q^{A_2}({_3M_2},{_2M_1}), {_3X_1})
\ \Rightarrow\ 
\mathbf{Ext}^{p+q}_{A_{3}-A_{1}}(\mathcal{C},{_3X_1}).$$
PREUVE. Soit $\mathcal{C}$ le complexe
$$\ldots\rightarrow {_3M_2}^q{_2M_1}\stackrel {d_q} \longrightarrow
{_3M_2}^{q-1}{_2M_1}\rightarrow\ldots\rightarrow
{_3M_2}^1{_2M_1}\rightarrow {_3M_2}\ {_2M_1}\rightarrow 0$$
de diff\'erentielle d\'efinie par
\begin{eqnarray*}
d_{q}(y,b_{1},\ldots
,b_{q},x)&=&-\,(yb_{1},b_{2},\ldots ,b_{q},x)\\
&\ +&\sum_{i=1}^{q-1}(-1)^{i+1}\,(y,b_{1},\ldots ,b_{i}b_{i+1},\ldots
,b_{q},x)\\
&\ +&(-1)^{q+1}(y,b_{1},\ldots
,b_{q-1},b_{q}x).
\end{eqnarray*}
L'homologie de $\mathcal C$ est alors 
$H_*(\mathcal{C})=Tor_*^{A_{2}}({_3M_2},{_2M_1}).$ 
En effet, $\mathcal{C}$ s'identifie au complexe
${_3M_2}\otimes_{A_{2}}Bar({_2M_1})$, où $Bar({_2M_1})$ est
la r\'esolution \emph{Bar} de $_2M_1$ comme $A_2$-module
\`a gauche. 
On consid\`ere maintenant le complexe double
$\mathcal{D}$ suivant au-dessus de $\mathcal{C}$:
$$\begin{array}{ccccccc}
 & & \vdots & & & & \vdots\\
 & & \downarrow & & & & \downarrow \\
\cdots & \rightarrow &
\bigoplus\limits_{\substack{{p+r=l+2}\\{p>0,\,r>0}}}
{^{r}{_3M_2}^{q}M_{1}^{p}} & \stackrel {d_{l,q}^{h}} \longrightarrow &
\cdots & \rightarrow
&\bigoplus\limits_{\substack{{p+r=l+2}\\{p>0,\,r>0}}}
{^{r}{_3M_2}\,{_2M_1}^{p}} \\
 & & \downarrow & & & & \downarrow \\
 & & \vdots & & & & \vdots\\
 & & \downarrow & & & & \downarrow \\
 \cdots & \rightarrow & {^{1}{_3M_2}^{q}{_2M_1}^{1}} &
 \stackrel {d_{0,q}^{h}} \longrightarrow &
\cdots & \rightarrow & {^1{_3M_2}\,{_2M_1}^1}\\
 & & \downarrow \epsilon_{q} & & & &
 \downarrow \epsilon_{0}\\
 \cdots & \rightarrow & {_3M_2}^{q}{_2M_1} & \stackrel
 {d_{q}} \longrightarrow &
\cdots & \rightarrow & {_3M_2}\,{_2M_1}\\
 & & \downarrow & & & & \downarrow \\
 & & 0 & & & & 0
\end{array}$$
La $q^e$ colonne de $\mathcal{D}$ est une résolution libre de
${_3M_2}^q{_2M_1}$, du m\^eme type que la résolution libre 
du lemme \ref{lem1}. L'expression des différentielles 
horizontales de $\mathcal D$ est semblable à celle de la 
différentielle de $\mathcal C$.

\begin{lem}\label{CE}
Le complexe double $\mathcal D$ est une
r\'esolution de Cartan-Eilenberg de $\mathcal C$ par des 
$A_3$-$A_1$-bimodules.
\end{lem}

Ce lemme suffit pour conclure la preuve de la
proposition: il reste \`a appliquer \`a $\mathcal D$ le foncteur
$Hom_{A_3-A_1}( -\ ,{_3X_1})$, et \`a prendre le
complexe total $Tot^{\oplus}\big(Hom_{A_3-A_1}
(\mathcal{D},{_3X_1})\big)$, qui est isomorphe à
la troisième colonne de la suite spectrale:
$$\Big(Tot^{\oplus}\big(Hom_{A_3-A_1}
(\mathcal{D},{_3X_1})\big)\Big)^* \cong Hom_{A_3-A_1}
\Big(\bigoplus\limits_{\substack{{p+q+r=*+2}\\{p>0,\,r>0}}}
 {^r{_3M_2}^q{_2M_1}^p}\,,\,{_3X_1}\Big)$$
 
$$\cong Hom_{R-R}
\Big(\bigoplus\limits_{p+q+r=*}
 {^r{_3M_2}^{q}{_2M_1}^p}\,,\,X\Big) 
 = E_0^{2,*}.$$
La cohomologie de ce complexe est alors $\mathbf{Ext}
^{*}_{A_3-A_1}(\mathcal C,{_3X_1})$.\\

PREUVE DU LEMME~\ref{CE}.
Il y a deux points \`a v\'erifier:
\begin{enumerate}
\item Pour tout $q\geq 0$, les bords $B^{h}_{*,q}(\mathcal{D})$ de la diff\'erentielle
horizontale de $\mathcal{D}$ forment une r\'esolution
projective de $B_{q}(\mathcal{C})$, le $q^{e}$ bord de $\mathcal{C}$.
\item Pour tout $q\geq 0$, les groupes d'homologie
$H^{h}_{*,q}(\mathcal{D})$ par rapport \`a la diff\'erentielle
horizontale de $\mathcal{D}$ forment une r\'esolution
projective de $H_{q}(\mathcal{C})$.
\end{enumerate}
Pour montrer
ces deux propriétés, il suffit de remarquer 
que la diff\'erentielle horizontale
$$d_{l,q}^{h}:\bigoplus_{\substack{{p+r=l+2}\\{p>0,\,r>0}}}
 {^{r}{_3M_2}^{q}{_2M_1}^{p}}\longrightarrow \bigoplus_{\substack{{p+r=l+2}\\{p>0,\,r>0}}}
 {^{r}{_3M_2}^{q-1}{_2M_1}^p}$$
 s'identifie \`a
$$\bigoplus_{\substack{{p+r=l+2}\\{p>0,\,r>0}}}
(-1)^{r+1}\,1_{A_{3}}\otimes\ldots\otimes 1_{A_{3}}\otimes
d_{q}\otimes 1_{A_{1}}\otimes\ldots\otimes 1_{A_{1}},$$
où $d_q$ est la $q^e$ différentielle de $\mathcal
C$. On a donc les identifications 
$$\left\{\begin{array}{l}
Z^h_{l,q}(\mathcal{D})=Ker\ d_{l,q}^{h}\ \cong\bigoplus\limits_
{\substack{{p+r=l+2}\\{p>0,\,r>0}}}
{^{r}\big(Ker\,d_{q}\big)^{p}}\ =\bigoplus\limits_
{\substack{{p+r=l+2}\\{p>0,\,r>0}}}
{^{r}\big(Z_q(\mathcal C)\big)^{p}}\\
\\
B^h_{l,q}(\mathcal{D})=Im\ d_{l,q+1}^{h}\ \cong\bigoplus\limits_
{\substack{{p+r=l+2}\\{p>0,\,r>0}}}
{^{r}\big(Im\,d_{q+1}\big)^{p}}\ =\bigoplus\limits_
{\substack{{p+r=l+2}\\{p>0,\,r>0}}}
{^{r}\big(B_q(\mathcal C)\big)^{p}}
\end{array}\right.$$
Ainsi, $B^{h}_{*,q}(\mathcal{D})\twoheadrightarrow 
B_{q}(\mathcal{C})$ est de la m\^eme forme que le complexe donné
au lemme \ref{lem1}; c'est donc une résolution libre de 
$B_{q}(\mathcal{C})$. On a de plus:
$$H^{h}_{l,q}(\mathcal{D})\ \cong\bigoplus\limits_
{\substack{{p+r=l+2}\\{p>0,\,r>0}}}
\frac{^r\big(Z_q(\mathcal C)\big)^p}
{^r\big(B_q(\mathcal C)\big)^p}
\ \cong\bigoplus\limits_
{\substack{{p+r=l+2}\\{p>0,\,r>0}}}{^{r}\Big(
\frac{Z_q(\mathcal C)}{B_q(\mathcal C)}\Big)^{p}}\ =\bigoplus
\limits_{\substack{{p+r=l+2}\\{p>0,\,r>0}}}
{^{r}\big(H_q(\mathcal C)\big)^{p}},$$
donc $H^{h}_{*,q}(\mathcal{D})\twoheadrightarrow 
H_{q}(\mathcal{C})$ est encore une r\'esolution 
libre. $\diamondsuit$\\

Nous avons donc complètement explicité les termes du 
niveau 1 de la suite spectrale; la ligne $l$ est
$$\begin{array}{l}
\quad HH^l(A_1,{_1X_1})\\
\oplus\ HH^l(A_2,{_2X_2})\\
\oplus\ HH^l(A_3,{_3X_3})
\end{array}
\stackrel{d_1^{0,l}}{\longrightarrow}
\begin{array}{l}
\quad Ext^l_{A_2-A_1}({_2M_1},{_2X_1})\\
\oplus\ Ext^l_{A_3-A_2}({_3M_2},{_3X_2})\\
\oplus\ Ext^l_{A_3-A_1}({_3M_1},{_3X_1})
\end{array}
\stackrel{d_1^{1,l}}{\longrightarrow}\mathbf{Ext}^l_{A_3-A_1}
(\mathcal{C},{_{3}X_{1}}).$$
Il nous reste maintenant à expliciter les deux 
différentielles. Pour cela, rappelons tout d'abord que pour 
$f\in Hom_k(A_1^{\otimes_kl},{_1X_1})$ et
$\theta \in End_{A_2-A_1}({_2M_1})$, on définit l'élément 
$\theta\smallsmile f$ de $Hom_k({_2M_1}
A_1^{\otimes_kl}, {_2X_1})$ de la manière suivante:
\begin{eqnarray*}
\theta\smallsmile f :\ {_2M_1}A_1^{\otimes_kl} & 
\longrightarrow & {_2X_1}\\
(x, a_1, \ldots , a_l) & \longmapsto & \theta(x)f(a_1,
\ldots , a_l),
\end{eqnarray*}
le produit ${_2M_1}\otimes{_1X_1}\longrightarrow{_2X_1}$
étant induit par la structure de $\mathcal T$-module à
gauche de $X$ $\mathcal T\otimes X\longrightarrow X$. 
Le produit \emph{cup} ainsi défini induit un produit en
cohomologie:
$$End_{A_2-A_1}({_2M_1})\otimes HH^l(A_1,{_1X_1})
\longrightarrow Ext_{A_2-A_1}^*({_2M_1},{_2X_1});$$
En effet, on vérifie facilement que $\delta(\theta\smallsmile
f)=-\theta\smallsmile\delta f$. De la m\^eme façon, 
pour $\theta\in End_{A_3-A_2}({_3M_2})$ et $f\in 
Hom_{A_2-A_1}({_2M_1}A_1^{\otimes_kl-1},{_2X_1})$, on peut 
construire un élément $\theta\smallsmile f$ de $Hom_{A_3-A_1}({_3M_2}\
{_2M_1}A_1^{\otimes_kl-1},{_3X_1})$ de la manière suivante:
\begin{eqnarray*}
\theta\smallsmile f :\ {_3M_2}\ {_2M_1}^{l-1} & 
\longrightarrow & {_3X_1}\\
(y,x,a_1,\ldots,a_{l-1}) & \longmapsto & \theta(y)
f(x,a_1,\ldots,a_{l-1});
\end{eqnarray*}
cette construction induit aussi un 
produit \emph{cup} en cohomologie:
$$End_{A_3-A_2}({_3M_2})\otimes Ext^l_{A_2-A_1}({_2M_1},{_2X_1})
\longrightarrow\mathbf{Ext}^l_{A_3-A_1}(\mathcal C,{_3X_1}).$$

\begin{prop}
La première diff\'erentielle du niveau 1\\
\begin{flushleft}
$d_{1}^{0,l}:HH^{l}(A_{1},{_{1}X_{1}})\oplus 
HH^{l}(A_{2},{_{2}X_{2}})\oplus 
HH^{l}(A_{3},{_{3}X_{3}})$
\end{flushleft}
\begin{flushright}
$\longrightarrow Ext^{l}_{A_{2}-A_{1}}({_2M_1},{_2X_1})\oplus
Ext^l_{A_3-A_2}({_3M_2},{_3X_2})\oplus
Ext^l_{A_3-A_1}({_3M_1},{_{3}X_{1}})$
\end{flushright}
est donn\'ee par la somme de produits \emph{cup}:
\begin{eqnarray*}
f+g+h & \longmapsto & \ \ \ 1_{_2M_1}\smallsmile f
+(-1)^{l+1}\, g\smallsmile 1_{_2M_1}\\
& & +\ 1_{_3M_2}\smallsmile g +(-1)^{l+1}\,
h\smallsmile 1_{{_3M_2}}\\
& & +\ 1_{_3M_1}\smallsmile f 
+(-1)^{l+1}\, h\smallsmile
1_{_3M_1}\ ,
\end{eqnarray*}
où $1_N$ désigne le morphisme identité du module $N$. 
La seconde diff\'erentielle
$$d_1^{1,l}:Ext^l_{A_2-A_1}({_2M_1},{_2X_1})\oplus
Ext^l_{A_3-A_2}({_3M_2},{_3X_2})\oplus
Ext^l_{A_3-A_1}({_3M_1},{_{3}X_{1}})\rightarrow \mathbf{Ext}^{l}_{A_3-A_1}
(\mathcal{C},{_3X_1})$$
est donn\'ee par 
\begin{eqnarray*}
f+g+h&\longmapsto&1_{{_3M_2}}\smallsmile f
+(-1)^{l+2}\, g\smallsmile 1_{{_2M_1}}\ +\ \delta h,
\end{eqnarray*} 
où $\delta h$ est d\'efinie pour $h:{^q{_3M_1}^p}
\longrightarrow{_3X_1}$ par $\delta h=(-1)^{q+1}
h\circ(1_{A_3}^{\otimes_k q}\otimes\mu\otimes
1_{A_1}^{\otimes_k p}):{^q{_3M_2}\ {_2M_1}^p}\longrightarrow
{_3X_1}$.
\end{prop}

Les résultats de cette proposition découlent de calculs
directs. Nous terminons l'étude du cas $n=3$ en donnant 
dans le cas
des algèbres triangulaires \emph{tensorielles} une condition
pour que la suite spectrale d\'eg\'en\`ere au niveau 2:

\begin{prop}\label{degenere}
La suite spectrale de cohomologie associée à l'algèbre 
triangulaire tensorielle
$$\left(\begin{array}{ccc}
A_1\\
_2M_1 & k\\
_3M_2\otimes\,{_2M_1} & _3M_2 & A_3
\end{array}\right)$$
dégénère au niveau 2.
\end{prop}

La d\'emonstration de cette proposition repose sur le lemme
suivant:

\begin{lem}\label{niveau 2}
Pour une algèbre triangulaire tensorielle de taille $3$
$$\left(\begin{array}{ccc}
A_1\\
_2M_1 & A_2\\
_3M_2\otimes_{A_2}{_2M_1} & _3M_2 & A_3
\end{array}\right)$$
la différentielle au niveau 2 de la suite spectrale de
cohomologie associée
$$d_2^{0,l}: E_2^{0,l}=Ker\ d_1^{0,l}\longrightarrow 
\frac{E_1^{2,l-1}}{Im\ d_1^{1,l}}=E_2^{2,l-1}$$
est nulle sur les composantes de $E_2^{0,l}$ d\'ependant de
$A_1$ et $A_3$:
$$d_2^{0,l}\Big(\big(HH^l(A_1,{_1X_1})\oplus
HH^l(A_{3},{_3X_3})\big)\cap Ker\ d_1^{0,l}\Big)=0.$$
\end{lem}

En effet, dans ce cas, la proposition \ref{degenere} en découle : 
comme $HH^*(k,{_2X_2})=0$ si $*>0$, le terme
$E_2^{0,*}$ pour $*>0$ est exactement $\big(HH^*(A_1,{_1X_1})
\oplus HH^*(A_{3},{_3X_3})\big)\cap Ker\ d_1^{0,*}$;
d'apr\`es le lemme la diff\'erentielle
$d_2^{0,*}$ est donc nulle.\\

PREUVE DU LEMME \ref{niveau 2}. 
Soit $f\in Hom_k(A_1^{\otimes l},{_1X_1})$ et $h\in
Hom_k(A_3^{\otimes l},{_3X_3})$ des
cocycles. La dif\-fé\-ren\-tiel\-le horizontale appliqu\'ee \`a
$f+h$ a quatre composantes: 
$$\left\{\begin{array}{cll}
1_{_2M_1}\smallsmile f &\in& Hom_k({_2M_1}A_1^{\otimes l}\ ,\ 
{_2X_1})\\
h\smallsmile 1_{_3M_2} &\in& Hom_k(A_3^{\otimes l}{_2M_1}\ ,\ 
{_3X_2})\\
1_{_3M_2\otimes_{A_2}{_2M_1}}\smallsmile f
&\in& Hom_k({_3M_2}\otimes_{A_2}{_2M_1}A_1^{\otimes l}\ ,\  
{_3X_1})\\
h\smallsmile 1_{_3M_2\otimes_{A_2}{_2M_1}}&\in& 
Hom_k(A_3^{\otimes l}{_3M_2}\otimes_{A_2}{_2M_1}\ ,\  {_3X_1}).
\end{array}\right.$$
Si $f+h$ est dans $Ker\ d_1^{0,l}$, alors ces quatre éléments 
sont des cobords par rapport à la différentielle verticale,
c'est-\`a-dire qu'il existe
$$\left\{\begin{array}{l}
\phi\ \in\ Hom_k({_2M_1} A_1^{\otimes l-1}\ ,\ {_{2}X_{1}})\\
\psi\ \in\ Hom_k(A_3^{\otimes l-1}{_3M_2}\ ,\ {_{3}X_{2}})\\
\chi\ \in\ Hom_k({_3M_2}\otimes_{A_{2}}{_2M_1}
A_1^{\otimes l-1}\ ,\ {_{3}X_{1}})\\
\xi\ \in\ Hom_k(A_3^{\otimes l-1}{_3M_2}\otimes_{A_{2}}{_2M_1}\ ,\ {_{3}X_{1}})\ ,
\end{array}\right.$$
tels que, pour $x\in {_2M_1}$, $y \in {_3M_2}$, 
$a_1,\ldots,a_l\in A_1$ et $c_1,\ldots,c_l\in A_3$, 
les \'egalit\'es suivantes soient v\'erifi\'ees:
$$\begin{array}{l}
\phi(xa_1,\ldots,a_l)+\sum\limits_{i=1}^{l-1}(-1)^i
\phi(x,\ldots,a_ia_{i+1},\ldots)+(-1)^l\phi(x,
\ldots,a_{l-1})a_l
= x f(a_1,\ldots,a_l)\\
c_1\psi(c_2,\ldots,y)+\sum\limits_{i=1}^{l-1}(-1)^i
\psi(\ldots,c_ic_{i+1},\ldots,y)
+(-1)^l\psi(c_1,\ldots,c_ly)=h(c_1,\ldots,c_l)y\\
\chi(y\!\otimes\!xa_1,\ldots,a_l)+\sum\limits_{i=1}^{l-1}(-1)^i
\chi(y\!\otimes\!x,\ldots,a_ia_{i+1},\ldots)
+(-1)^l\chi(y\!\otimes\!x,a_1,\ldots)a_l
=y\!\otimes\!xf(a_1,\ldots,a_l)\\
c_1\xi(c_2,\ldots,c_l,y\!\otimes\!x)+\sum\limits_{i=1}^{l-1}(-1)^i
\xi(\ldots,c_ic_{i+1},\ldots,y\!\otimes\!x)
+(-1)^l\xi(c_1,\ldots,c_ly\!\otimes\!x)=h(c_1,\ldots,c_l)y\!\otimes\!x
.\end{array}$$
Pour obtenir un repr\'esentant de $d_2^{0,l}(f+h)$, il reste 
\`a appliquer \`a $\phi$, $\psi$, $\chi$ et $\xi$ la
deuxi\`eme diff\'erentielle horizontale; on obtient ainsi que
$$d_2^{0,l}(f+h)\ :\ {_3M_2}\ {_2M_1}A_1^{l-1}\oplus
A_3^{l-1}{_3M_2}\ {_2M_1}\longrightarrow{_{3}X_{1}}$$
est donné par:
$$\begin{array}{l}
d_2^{0,l}(f+h)\big((y,x,a_1,\ldots,a_{l-1})+(c_1,\ldots,c_{l-1},y',x')\big)\\
\\
=y\phi(x,a_1,\ldots,a_{l-1})-\chi(y\!\otimes\!x,a_1,\ldots,a_{l-1})
+(-1)^l\big(\psi(c_1,\ldots,c_{l-1},y')x'-\xi(c_1,\ldots,c_{l-1},y'\!\otimes\!x')\big).
\end{array}$$
Or $d_2^{0,l}(f+h)$ ne dépend pas du choix des quatre 
\'el\'ements $\phi$, $\psi$, $\chi$ et $\xi$. 
Remarquons alors que l'on peut définir
$\chi$ et $\xi$ en fonction respectivement de $\phi$ 
et de $\psi$, de fa\c con à ce qu'ils vérifient 
les égalités demandées, en posant:
$$\left\{\begin{array}{lll}
\chi(y\!\otimes\!x,a_1,\ldots,a_{l-1})&=&
y\phi(x,a_1,\ldots,a_{l-1})\\
\xi(c_1,\ldots,c_{l-1},y\!\otimes\!x)&=&
\psi(c_1,\ldots,c_{l-1},y)x\ .
\end{array}\right.$$
Les morphismes $\chi$ et $\xi$ sont bien 
d\'efinis, car $\phi$ et $\psi$ sont $A_2$-lin\'eaires,
respectivement  à gauche et à droite. En effet, la 
différentielle verticale appliquée à $\phi$ a
deux composantes: l'une dans $Hom_k({_2M_1} A_1^{\otimes l},
{_2X_1})$, qui est par hypothèse $1_{_2M_1}\smallsmile f$, et
l'autre dans $Hom_k(A_2\ {_2M_1}A_1^{\otimes l-1},{_2X_1})$,
qui doit \^etre nulle. Cela signifie explicitement que
l'application
\begin{eqnarray*}
A_2\ {_2M_1}A_1^{l-1}&\longrightarrow&{_{2}X_{1}}\\
(b,x,a_1,\ldots,a_{l-1})&\longmapsto&
b\phi(x,a_1,\ldots,a_{l-1}) - \phi(bx,a_1,\ldots,a_{l-1})
\end{eqnarray*}
est nulle; ainsi $\phi$ est bien $A_2$-linéaire à gauche. On
montre de la m\^eme fa\c con que $\psi$ est $A_2$-linéaire à
droite. 
En écrivant alors $\chi$ et $\xi$ sous la forme ci-dessus dans
l'expression de $d_2^{0,l}(f+h)$, on
obtient $d_2^{0,l}(f+h)=0$. $\diamondsuit$

\subsection{La suite spectrale pour $n$ quelconque}

Nous allons terminer cette section en donnant des
r\'esultats pour $n\geq 3$. Sous des hypothèses
de projectivité sur certains modules, les termes du niveau 1
de la suite spectrale seront des sommes de groupes
d'extension. Dans le cas général, il suffirait de remplacer
dans les calculs les foncteurs $Ext$ par des foncteurs $HyperExt$ embo\^\i tés.
Par ailleurs, les différentielles au niveau 1 s'expriment
encore à l'aide de produits \emph{cup}.

Rappelons que les termes au niveau 0 de la suite
spectrale sont 
$$E_0^{t,l-t}=Hom_{R-R}\Bigg( \bigoplus_{\substack{\tau\in
TR_{l}(\mathcal{Q}) \\ |\tau|=t}} M_{\tau}\ ,\ X\Bigg),$$
où le module $M_{\tau}$, associé à la
$l$-trajectoire $\tau$ de longueur $t$, est de la
forme
$$M_{\tau}={^{p_{t+1}}{_{k_{t+1}}M_{k_t}}^{p_{t}}\ldots
 ^{p_{2}}{_{k_2}M_{k_1}}^{p_{1}}}\ ,\ \textrm{avec}
 p_{1}+\ldots +p_{t+1}=l-t.$$
Supposons alors que pour tout i variant de $2$ à
$t$, pour tous entiers $p_{i+1},\ldots,p_t$, les 
$A_{k_{i}}^{e}$-modules
$$_{k_{t+1}}M_{k_t}^{p_{t}}\ldots
^{p_{i+1}}{_{k_{i+1}}M_{k_i}}\ 
{_{k_i}M_{k_{i-1}}}\otimes_{A_{k_{i-1}}}\ldots\otimes_{A_{k_2}}
{_{k_2}M_{k_1}}$$
sont projectifs.

\begin{prop}\label{ext}
Sous l'hypothèse de projectivité ci-dessus, le
niveau 1 de la suite spectrale est donné par
$$\left\{\begin{array}{ccl}
E_{1}^{0,\,*}&=&HH^{*}(A,X_{0})\\
\\
E_{1}^{1,\,*}&=&\bigoplus\limits_{\substack{{1\leq i\leq
n-1}\\{1\leq \alpha\leq n-i}}}
\quad Ext^{*}_{A_{i+\alpha}-A_i}({_{i+\alpha}M_i}\ ,\ {_{i+\alpha}X_i})\\
\\
E_{1}^{t,\,*}&=& \bigoplus\limits_
{1\leq k_1<\ldots<k_{t+1}\leq n}
Ext^{*}_{A_{k_{t+1}}-A_{k_1}}\ \big({_{k_{t+1}}M_{k_t}}\otimes_
{A_{k_t}}\ldots\otimes_{A_{k_2}}{_{k_2}M_{k_1}}\,,\,{_{k_{t+1}}X_{k_1}}
\big),\ t\geq 2.\end{array}\right.$$
\end{prop}

Le calcul des termes $E_{1}^{0,\,*}$ et $E_{1}^{1,\,*}$ est 
le m\^eme que dans le cas $n=3$ et ne demande pas
d'hypoth\`eses suppl\'ementaires. En revanche, nous
allons utiliser le lemme suivant pour le calcul de 
$E_{1}^{t,\,*}$ pour $t\geq2$:

\begin{lem}\label{lem2}
Soient $A$ une $k$-alg\`ebre, $N$ et $N'$ des $A$-modules
 respectivement \`a droite et \`a gauche. Si
 \emph{$N\!\otimes_{k}\!N'$ est projectif comme
$A^{e}$-module}, alors le complexe suivant:\\

\begin{tabular}{lcc}
$\mathcal{R}$: & $\qquad\qquad$ & $\ldots\rightarrow N^{l+1}N' \stackrel {d_{l}} 
{\longrightarrow} N^{l}N'\rightarrow
\ldots\rightarrow N^{1}N'\rightarrow
N\,N'\rightarrow 0$
\end{tabular}\\
\\
de diff\'erentielle
\begin{eqnarray*}
d_{l}:\qquad {N^{l+1}N'} & \longrightarrow &
{N^{l}N'}\\
(x,a_{1},\ldots,a_{l},x') & \longmapsto
& \quad (xa_{1},a_{2},\ldots,a_{l},x')\\
& & +\ \sum_{i=1}^{l-1}(-1)^{i}\
(x,a_{1},\ldots,a_{i}a_{i+1},\ldots,a_{l},x')\\
& & +\ (-1)^{l}\ (x,a_{1},\ldots,a_{l-1},a_{l}x')
\end{eqnarray*}
est acyclique en degr\'e $\geq0$, et son homologie
en degr\'e $-1$ est $N\!\otimes_{A}\!N'$.
\end{lem}

PREUVE.
Il suffit de remarquer que $\mathcal{R}$ s'identifie
au complexe $(N\otimes_{k}N')\otimes_{A^{e}}Bar(A)$, de 
diff\'erentielle $d_*=-1_{N\otimes_kN'}\otimes\partial_{*+1}$,
où $(Bar,\partial)$ est la résolution $Bar$ de $A$. 
Le $A^{e}$-module $N\otimes_{k}N'$ étant suppos\'e projectif, 
le foncteur $(N\otimes_{k}N')\otimes_{A^{e}}-$ est exact. 
Comme $H_{l}(Bar_{*}(A))=0$ pour $l>0$ et
$H_{0}(Bar_{*}(A))=A$, on obtient:
$$\left\{ \begin{array}{l}
H_{l}(\mathcal{R})=0 \quad pour\ l\geq 0,\\
H_{-1}(\mathcal{R})=N\!\otimes_{k}\!N'\!\otimes_{A^{e}}\!A\cong
N\!\otimes_{A}\!A\otimes_{A}\!N'\cong
N\!\otimes_{A}\!N'.\ \diamondsuit
\end{array}\right.$$
PREUVE DE LA PROPOSITION~\ref{ext}. Soit une suite d'entiers 
fixée $1\leq k_1<\ldots<k_{t+1}\leq n$. 
On reprend le complexe du lemme~\ref{lem1}, en
rempla\c cant $N$ par ${_{k_{t+1}}M_{k_t}}^{p_{t}}\ldots
{^{p_3}{_{k_3}M_{k_2}}^{p_{2}}{_{k_2}M_{k_1}}}$:
$$\ldots\rightarrow\bigoplus_{\substack{{p_{t+1}+p_{1}=l+2}\\{p_{t+1}>0,\,
p_{1}>0}}}
{^{p_{t+1}}{_{k_{t+1}}M_{k_t}}^{p_{t}}\ldots
 ^{p_{2}}{_{k_2}M_{k_1}}^{p_{1}}} 
\rightarrow\cdots
\rightarrow{_{k_{t+1}}M_{k_t}}^{p_{t}}\ldots
 ^{p_{2}}{_{k_2}M_{k_1}}\rightarrow 0.$$
C'est une r\'esolution libre de
${_{k_{t+1}}M_{k_t}}^{p_{t}}\ldots{^{p_3}{_{k_3}M_{k_2}}
^{p_{2}}{_{k_2}M_{k_1}}}$. 
On fait maintenant varier $p_{2}$ dans le complexe 
double suivant:
$$\begin{array}{cccccccccc}
 & & & \vdots & & & & \vdots & & \\
 & & & \downarrow & & & & \downarrow & & \\
\cdots\!\! & \!\rightarrow\! & \bigoplus\!\!\!\! &
^{p_{t+1}}{_{k_{t+1}}M_{k_t}}^{p_{t}}\ldots
 ^{p_{2}}{_{k_2}M_{k_1}}^{p_{1}}\!
& \!\rightarrow\! & \!\cdots\! &
\!\!\rightarrow\!\! & \!{_{k_{t+1}}M_{k_t}}^{p_{t}}\ldots
 ^{p_{2}}{_{k_2}M_{k_1}}\! & \!\rightarrow\! 0\\
 & & {\substack{{p_{t+1}+p_{1}=l+2}\\{p_{t+1}>0,\,
p_{1}>0}}} & & & & & & & \\
 & & & \downarrow & & & & \downarrow & & \\
 & & & \vdots & & & & \vdots & & \\
 & & & \downarrow & & & & \downarrow & & \\
\cdots\!\! & \!\rightarrow\! & \bigoplus\!\!\!\! &
{^{p_{t+1}}{_{k_{t+1}}M_{k_t}}^{p_{t}}\ldots
 ^1{_{k_2}M_{k_1}}^{p_{1}}}\! 
& \!\rightarrow\! & \!\cdots\! &
\!\!\rightarrow\!\! &
\!{_{k_{t+1}}M_{k_t}}^{p_{t}}\ldots
 ^1{_{k_2}M_{k_1}}
\! & \!\rightarrow\! 0\\
  & & {\substack{{p_{t+1}+p_{1}=l+2}\\{p_{t+1}>0,\,
p_{1}>0}}} & & & & & & & \\
 & & & \downarrow & & & & \downarrow & &  \\
\cdots\!\! & \!\rightarrow\! & \bigoplus\!\!\!\! &
{^{p_{t+1}}{_{k_{t+1}}M_{k_t}}^{p_{t}}\ldots
{_{k_2}M_{k_1}}^{p_{1}}} \!
& \!\rightarrow\! & \!\cdots\! &
\!\!\rightarrow\!\! & \!{_{k_{t+1}}M_{k_t}}^{p_{t}}\ldots
{_{k_2}M_{k_1}}\! & \!\rightarrow\! 0\\
& & {\substack{{p_{t+1}+p_{1}=l+2}\\{p_{t+1}>0,\,
p_{1}>0}}} & & & & & & & \\
  &  & & \downarrow &  &  &  & \downarrow &  &  \\
  &  & & 0 &  &  & & 0 &  &
\end{array}$$
Les lignes de ce complexe double sont des résolutions libres
du type de celle donnée au lemme \ref{lem1}. En
particulier, 
les lignes sont acycliques, donc le complexe total
$\mathcal{B}_{2}$ l'est aussi.
On consid\`ere alors le sous-complexe
$\mathcal{S}_{2}$ form\'e par la derni\`ere
colonne. On lui applique le lemme~\ref{lem2} en
supposant $_{k_{t+1}}M_{k_t}^{p_{t}}\ldots
^{p_{3}}{_{k_{3}}M_{k_2}}\ {_{k_2}M_{k_1}}$
projectif comme $A_{k_{2}}^{e}$-module. Ainsi l'homologie 
de $\mathcal S_2$ est nulle en degr\'e $\geq 0$, et vaut 
$_{k_{t+1}}M_{k_t}^{p_{t}}\ldots
^{p_{3}}{_{k_{3}}M_{k_2}}\otimes_{A_{k_2}}{_{k_2}M_{k_1}}$ 
en degr\'e $-1$. 
De la suite exacte longue d'homologie associ\'ee \`a
la suite exacte courte de complexes
$0\rightarrow\mathcal S_2\rightarrow\mathcal B_2\rightarrow
\mathcal B_2/\mathcal S_2\rightarrow 0$, on tire:
$$\left\{ \begin{array}{l}
H_{l}(\mathcal{B}_{2}/\mathcal{S}_{2})=0\quad
pour\ l\geq 0,\\
H_{-1}(\mathcal{B}_{2}/\mathcal{S}_{2})=_{k_{t+1}}M_{k_t}^{p_{t}}\ldots
^{p_{3}}{_{k_{3}}M_{k_2}}\otimes_{A_{k_2}}{_{k_2}M_{k_1}}.
\end{array}\right.$$
On augmente alors $\mathcal B_2/\mathcal S_2$ par $_{k_{t+1}}M_{k_t}^{p_{t}}\ldots
^{p_{3}}{_{k_{3}}M_{k_2}}\otimes_{A_{k_2}}{_{k_2}M_{k_1}}$,
et on d\'ecale les indices d'un cran. On a donc obtenu une
résolution libre de $_{k_{t+1}}M_{k_t}^{p_{t}}\ldots
^{p_{3}}{_{k_{3}}M_{k_2}}\otimes_{A_{k_2}}{_{k_2}M_{k_1}}$:
$$\ldots\rightarrow\bigoplus_{\substack{
{p_{t+1}+p_2+p_1=l+2}\\{p_{t+1}>0,\,p_1>0}}}\!\!\!
{^{p_{t+1}}}_{k_{t+1}}M_{k_t}^{p_{t}}\ldots
^{p_2}{_{k_2}M_{k_1}}^{p_{1}}\rightarrow\ldots
\rightarrow {_{k_{t+1}}M_{k_t}^{p_{t}}\ldots
^{p_{3}}{_{k_{3}}M_{k_2}}\otimes_{A_{k_2}}{_{k_2}M_{k_1}}}\rightarrow 0.$$
On it\`ere alors ce proc\'ed\'e, en faisant varier
successivement $p_{i}$ ($2\leq i\leq t$) dans un complexe double
$\mathcal{B}_{i}$, en supposant que 
$_{k_{t+1}}M_{k_t}^{p_{t}}\ldots
^{p_{i+1}}{_{k_{i+1}}M_{k_i}}\ 
{_{k_i}M_{k_{i-1}}}\otimes_{A_{k_{i-1}}}\ldots\otimes_{A_{k_2}}
{_{k_2}M_{k_1}}$
est projectif comme $A_{k_{i}}^{e}$-module. On 
obtient finalement une r\'esolution libre de  
$_{k_{t+1}}M_{k_t}\otimes_{A_{k_t}}\ldots\otimes_{A_{k_2}}
{_{k_2}M_{k_1}}$ comme $A_{k_{t+1}}$-$A_
{k_{1}}$-bimodule:
$$\ldots\rightarrow\bigoplus_{\substack{{p_{t+1}+\ldots+p_{1}=l+2}\\
{p_{t+1}>0,\,p_{1}>0}}}{^{p_{t+1}}{_{k_{t+1}}M_{k_t}}^{p_t}\ldots
{^{p_2}_{k_2}M_{k_1}}^{p_{1}}}\rightarrow\ldots\rightarrow
{_{k_{t+1}}M_{k_t}}\otimes_{A_{k_t}}\ldots\otimes_{A_{k_2}}
{_{k_2}M_{k_1}}\rightarrow 0.$$
Pour finir, on applique \`a cette r\'esolution le foncteur
$Hom_{A_{k_{t+1}}-A_{k_{1}}}(\ -\ ,
\,{_{k_{t+1}}X_{k_{1}}})$ et on simplifie en
utilisant la formule d'adjonction habituelle. On retrouve
ainsi une des composantes du terme $E_0^{t,l}$ de la suite spectrale, dont la
cohomologie est donc $$Ext^*_{A_{k_{t+1}}-A_{k_{1}}}\Big({_{k_{t+1}}M_{k_t}}\ldots
{_{k_2}M_{k_1}} ,
\,{_{k_{t+1}}X_{k_{1}}}\Big).\ \diamondsuit$$

Pour terminer, nous explicitons les
diff\'erentielles au niveau 1 de la suite spectrale.

\begin{prop}
La premi\`ere diff\'erentielle au niveau 1 est:
\begin{eqnarray*}
d_{1}^{0,l}\!:\bigoplus_{1\leq j \leq n}\!\!
HH^{l}(A_{j},{_{j}X_{j}})\!\!\! &
\!\!\rightarrow\!\!\!\!\!\! & 
\!\!\!\!\bigoplus_{\substack{{1\leq i\leq
n-1}\\{1\leq\alpha\leq n-i}}}\!
Ext^{l}_{A_{i+\alpha}-A_{i}}
\big({_{i+\alpha}M_i}\,,\,{_{i+\alpha}X_{i}}\big)\\
f_{1}+\ldots +f_{n} & \mapsto & 
\sum_{j=1}^{n}\ 
\bigg\{\ \sum_{\beta=1}^{n-j}\ 1_{_{j+\beta}M_j}\smallsmile f_{j}
+(-1)^{l+1}\sum_{\beta =1}^{j-1}\ f_{j}\smallsmile
1_{_jM_{j-\beta}}\ \bigg\}
\end{eqnarray*}
\end{prop}

Nous allons maintenant donner une expression des
diff\'erentielles $d_{1}^{t,l}$ pour $t>0$; pour
cela, rappelons que l'on note pour $1\leq i<j<l\leq n$
$$\mu_{l,j,i}:\ {_lM_j}\otimes{_jM_i}\longrightarrow
{_lM_i}$$ les morphismes qui permettent de définir le produit 
de $\mathcal T$.

\begin{prop}
Soit
$$f:\ {^{p_{t+1}}{_{k_{t+1}}M_{k_t}}^{p_t}}\ldots
{^{p_2}{_{k_2}M_{k_1}}^{p_1}}\longrightarrow
{_{k_{t+1}}X_{k_1}}\ ,\ \textrm{ avec }p_1+\ldots+p_{t+1}=l-t,$$
un représentant d'un élément de la composante
$Ext_{A_{k_{t+1}}-A_{k_t}}^{l}\big({_{k_{t+1}}M_{k_t}}\otimes_{A_{k_t}}
\ldots\otimes{A_{k_2}M_{k_1}},
{_{i+\alpha}X_i}\big)$ de $E_1^{t,l}$. La différentielle 
$d_1^{t,l}$, restreinte à cette composante, est alors 
donnée par:
\begin{eqnarray*}
d_1^{t,l}f&=&\sum_{k_{t+2}=k_{t+1}+1}^{n}
1_{_{k_{t+2}}M_{k_{t+1}}}\smallsmile f
\ +\ \sum_{k_0=1}^{k_1-1}(-1)^{l+2} f\smallsmile
1_{_{k_1}M_{k_0}}\\
& +& \sum_{i=1}^t
\sum_{\alpha=k_i+1}^{k_{i+1}-1}(-1)^{p_{t+1}+\ldots +p_{i+1}+t-i+1}
f\circ(1^{\otimes p_{t+1}+\ldots p_{i+1}+t-i}\otimes
\mu_{k_{i+1},\alpha,k_i}\otimes 1^{\otimes p_i+\ldots +p_1+i-1}).
\end{eqnarray*}
\end{prop}

\section{Applications}\label{application}

\subsection{Algèbre des chemins d'un carquois sans cycles
orientés}

Nous considérons dans cette section l'algèbre sur un corps $k$ 
associée à un carquois sans cycles orientés. Nous montrons
qu'elle est triangulaire, et que sa cohomologie peut \^etre
calculée en utilisant la suite spectrale
construite à la section \ref{suite spectrale}.

\begin{defi}
Soit $\mathcal Q$ un carquois, et $k$ un corps. L'\emph{algèbre des 
chemins de $\mathcal Q$ sur $k$}, notée $k\mathcal Q$, est la
$k$-algèbre de base l'ensemble des chemins de $\mathcal Q$,
et dont le produit est donné sur la base par la composition 
des chemins quand elle est possible, et par 0 sinon.
\end{defi}

Nous allons considérer dans la suite deux types de carquois:
les \emph{carquois sans cycles orientés} et les 
\emph{carquois à niveaux}.
Un cycle orienté dans un carquois est un chemin dont le but
co\"\i ncide avec la source. Un carquois sera dit à $n$ 
niveaux si on peut indexer ses sommets par
$e_1^1,\ldots,e_1^{i_1},\ldots,e_n^1,\ldots,e_n^{i_n}$, de
fa\c con à ce qu'il n'existe pas de flèche entre les sommets
$e_r^i$ et $e_s^j$ si $r>s$.

\begin{prop}
Soit $\mathcal Q$ un carquois connexe sans cycles orientés.
Alors l'algèbre des chemins de $\mathcal Q$ est isomorphe à
une algèbre triangulaire.
\end{prop}

Cette proposition est une conséquence immédiate des deux 
lemmes suivants:

\begin{lem}\label{sanscycle}
Si $\mathcal Q$ est un carquois connexe sans cycles
orientés, alors on peut ordonner les sommets de $\mathcal Q$
de fa\c con à obtenir un carquois à niveaux.
\end{lem}

\begin{lem}\label{niv}
Si $\mathcal Q$ est un carquois à niveaux, alors l'algèbre
$k\mathcal Q$ est triangulaire.
\end{lem}

PREUVE DU LEMME \ref{sanscycle}
Rappelons tout d'abord qu'un carquois est un graphe orienté
\emph{fini}; soit $N$ le nombre de sommets du carquois sans cycles
orientés $\mathcal Q$. Nous allons montrer que l'on peut munir
$\mathcal Q$ d'une structure à $N$ niveaux; cette construction 
n'est pas canonique, et, dans la pratique, le nombre de niveaux 
de $\mathcal Q$ pourra \^etre inférieur. Remarquons qu'il existe 
(au moins) un sommet de 
$\mathcal Q$ qui n'est le but d'aucune flèche; nous 
dirons qu'un tel sommet est \emph{initial}. Une construction des
niveaux de $\mathcal Q$ se fait par induction sur les
sommets initiaux. Soit $e$ un sommet initial; on définit le
niveau 1 de $\mathcal Q$ comme le niveau de $e$. On considère
alors le carquois dont l'ensemble des sommets est $\mathcal
Q\setminus \{e\}$. Ce carquois n'a pas de cycles orientés, il
possède donc un sommet initial $e'$. On définit le niveau 2 de
$\mathcal Q$ comme le niveau de $e'$. On itère ce procédé
jusqu'au $N^e$ sommet de $\mathcal Q$, qui formera le niveau
$N$. $\diamondsuit$\\

PREUVE DU LEMME \ref{niv}
Si $\mathcal Q$ a $n$ niveaux, on note pour tout
$r\in\{1,\ldots,n\}$ $e_r^1,\ldots,e_r^{i_r}$ les sommets 
formant le niveau $r$ de $\mathcal Q$, et on considère la
$k$-algèbre $A_r$ dont la base comme $k$-espace
vectoriel est $\{e_r^1,\ldots,e_r^{i_r}\}$.
De plus, pour $r\in\{1,\ldots,n-1\}$, et
$s\in\{r+1,\ldots,n\}$, on considère 
$_sM_r$ le $k$-espace vectoriel dont la base est
l'ensemble des chemins de tous les $e_r^i$ vers
les $e_s^j$. Alors l'algèbre $k\mathcal Q$ 
est l'algèbre triangulaire
$$\left(\begin{array}{cccc}
A_1\\_2M_1&A_2\\\vdots&\ddots&\ddots\\
_nM_1&\cdots&_nM_{n-1}&A_n\end{array}\right),$$
les applications
$\mu_{}:{_tM_s}\otimes_k{_sM_r}\rightarrow{_tM_r}$ étant
induites par la composition des chemins de $\mathcal Q$ 
quand elle est possible. $\diamondsuit$\\

Nous allons maintenant expliciter le calcul de la suite
spectrale sur un exemple de carquois sans cycles orientés:
\begin{diagram}[h=2em,w=3.5em]
&&c\\
&&\uTo\uTo\\
\mathcal Q:\qquad a&\rTo&b\\
&&\dTo\dTo\\
&&d
\end{diagram}	
L'algèbre $k\mathcal Q$ des chemins de $\mathcal Q$ est 
triangulaire tensorielle de taille 3:
$$k\mathcal Q=\left(\begin{array}{llll}
k\\
k & k\\
k^4 & k^4 & k^2
\end{array}\right).$$
La cohomologie de Hochschild de
$k\mathcal Q$ à valeurs dans $k\mathcal Q$ peut \^etre
calculée à l'aide de la 
suite spectrale construite à la
section \ref{suite spectrale}. 
Soit $\mathcal H$ le complexe de Hochschild relatif de
$k\mathcal Q$ par rapport à la sous-algèbre séparable 
$$R=\left(\begin{array}{llll}
k\\
&k\\
&& k^2
\end{array}\right),$$
dont le terme de degré $l$ est $\mathcal
H^{l}=Hom_{R-R}\big((k\mathcal Q)^{\otimes_{R}l},
k\mathcal Q\big)$, et 
dont la cohomologie est $HH^{*}(k\mathcal Q)$.
La suite spectrale associée à la filtration de $\mathcal H$
par la longueur des trajectoires de $\mathcal Q$ n'a au
niveau 1 qu'une seule ligne non nulle:
$$\begin{array}{ccccc}
R=k^4 &\stackrel {d^{0}} \longrightarrow & k\oplus End_{k^{2}-k}(k^{4})\oplus
End_{k^{2}-k}(k^4)
& \stackrel {d^{1}} \longrightarrow & End_{k^{2}-k}(k^4)\ .
\end{array}$$
En effet, les groupes d'extension de puissances de $k$ sont
nuls en degré positifs.
Les deux diff\'erentielles sont:
$$\begin{array}{rcl}
d^{0}:\qquad k^{4} & \longrightarrow & k\oplus End_{k^{2}-k}(k^{4})\oplus
End_{k^{2}-k}(k^4)\\
(\alpha, \beta, \gamma, \delta)
& \longmapsto & 
\Big((\alpha-\beta)\,id_{k^4},\,\phi,\,\psi\Big)\\
\end{array}$$
où $\phi$ et $\psi$ sont les applications:
$$\phi=\left(\begin{array}{cc}
(\beta-\gamma)id_{k^2}&0\\
0&(\beta-\delta)id_{k^2}
\end{array}
\right),\ 
\psi=\left(\begin{array}{cc}
(\alpha-\gamma)id_{k^2}&0\\
0&(\alpha-\delta)id_{k^2}
\end{array}
\right),$$
et
$$\begin{array}{rcl}
d^{1}:\quad k\oplus End_{k^{2}-k}(k^{4})\oplus
End_{k^{2}-k}(k^4) & \longrightarrow & End_{k^{2}-k}(k^4)\\
(\lambda,\phi,\psi) & \longmapsto & \lambda\,id_{k^4}
\ +\ \phi\ -\ \psi.
\end{array}$$
On a ainsi:
$$\left\{\begin{array}{ll}
Ker\ d^{0}\cong k\ ,&Ker\ d^{1}=k\oplus End_{k^{2}-k}(k^{4})\ ,\\
Im\ d^{0}\cong k^{3}\ ,&Im\ d^{1}=End_{k^{2}-k}(k^{4}).
\end{array}\right.$$
D'o\`u les termes de niveau 2 de la suite spectrale, et donc les
groupes $HH^{*}(k\mathcal Q)$:
$$\begin{array}{lll}
HH^{0}(k\mathcal Q)\cong k\ , & HH^{1}(k\mathcal Q)
\cong k^6,\ & HH^{*}(k\mathcal Q)=0
\text{ pour }*\geq 2.
\end{array}$$

\subsection{Liens avec des résultats antérieurs}

La suite spectrale construite à la section \ref{suite
spectrale}, et l'application aux carquois sans cycles
orientés, généralisent un résultat de Cibils sur l'algèbre
d'incidence d'un complexe simplicial. 
Soit $\Sigma$ un complexe simplicial de dimension finie.
On associe à $\Sigma$ un carquois $\mathcal Q_\Sigma$: les 
sommets sont 
les simplexes de $\Sigma$, et il existe une flèche entre le
sommet $\sigma$ et le sommet $\tau$ si et seulement si $dim\
\tau=\dim\ \sigma +1$ et le simplexe $\sigma$
est contenu dans le simplexe $\tau$. Ce carquois a des
niveaux naturels donnés par la dimension des simplexes de
$\Sigma$, ainsi l'algèbre $k\mathcal Q_\Sigma$ est 
triangulaire. Soit $I$ l'idéal de $k\mathcal Q_\Sigma$ 
engendré par les différences $\gamma-\gamma'$, où $\gamma$
et $\gamma'$ sont des chemins de m\^emes extrémités (de tels 
chemins sont dits \emph{parallèles}). L'algèbre $k\mathcal
Q_\Sigma/I$ est alors l'algèbre d'incidence du poset associé
à $\Sigma$ sur $k$; elle est encore triangulaire. D'après 
Gerstenhaber et Schack \cite{GS3}, la cohomologie
de Hochschild de l'algèbre $k\mathcal Q_\Sigma/I$ à valeurs dans
elle-m\^eme est isomorphe à la cohomologie simpliciale de
$\Sigma$ à coefficients dans $k$. Utilisant cet
isomorphisme, Cibils a construit dans \cite{Cib2} une suite
spectrale de cohomologie convergeant vers $H^*(\Sigma,k)$,
dont le niveau 1 est concentré sur la première ligne. 
Ce résultat peut aussi se déduire de la suite spectrale 
construite à la section \ref{suite spectrale}; en effet, comme
dans l'exemple précédent, la nullité des groupes d'extension
supérieurs de puissances de $k$ implique que seule la première
ligne du niveau 1 de la suite spectrale est non nulle.\\

Par ailleurs, si $\mathcal C$ est une catégorie finie, on
peut lui associer la $k$-algèbre $k\mathcal C$ des
morphismes de $\mathcal C$: sa base comme $k$-espace vectoriel
est l'ensemble des morphismes de $\mathcal C$, et son produit
est induit par la composition des morphismes. Par
exemple, si l'on considère la catégorie $\mathcal C_{\mathcal Q}$ induite par 
un carquois $\mathcal Q$, les objets de $\mathcal C_{\mathcal Q}$ étant les 
sommets de $\mathcal Q$, et les morphismes de $\mathcal C_{\mathcal Q}$ les 
chemins dans $\mathcal Q$, alors l'algèbre des morphismes de
$\mathcal C_{\mathcal Q}$ co\i\" ncide avec l'algèbre des chemins de
$\mathcal Q$. On peut définir la
cohomologie de Hochschild d'une catégorie $\mathcal C$ comme
la cohomologie de Hochschild de l'algèbre $k\mathcal C$. Bendiffalah et 
Guin ont étudié dans \cite{BG} la cohomologie des 
catégories \emph{musclées}. La suite spectrale construite à
la section \ref{suite spectrale} permet de calculer la
cohomologie des catégories dont l'algèbre des morphismes est
triangulaire. Ces deux 
situations co\"\i ncident pour une catégorie \emph{muscle}, 
dont l'algèbre est triangulaire de taille 2.

\end{document}